\documentclass[runningheads]{llncs}
\usepackage[T1]{fontenc}
\usepackage{graphicx}
\usepackage{lineno}
\usepackage{multirow}
\usepackage{adjustbox,amsmath,amssymb,upgreek}
\usepackage{bm}
\usepackage{pifont}
\usepackage[english]{babel}
\usepackage{geometry}
\usepackage[colorlinks,citecolor=blue]{hyperref}

\usepackage[
    backend=biber,
    style=nature,
  ]{biblatex}

\begin{document}
\title{Pupil-driven quantitative differential phase contrast imaging}
%
%
\author{Shuhe Zhang\inst{1,2,*} \and
Hao Wu\inst{1} \and
Tao Peng\inst{1} \and
Zeyu Ke\inst{1} \and
Meng Shao\inst{1} \and
Tos T. J. M. Berendschot\inst{2} \and
Jinhua Zhou\inst{1,3,**}}
\authorrunning{ }
%
\institute{School of Biomedical Engineering, Anhui Medical University, Hefei 230032, China \and
University Eye Clinic Maastricht, Maastricht University Medical Center + , P.O. Box 5800, Maastricht, 6202 AZ, the Netherlands\\
\and
Anhui Provincial Institute of Translational Medicine, Anhui Medical University, Hefei 230032, China\\
\email{* shuhe.zhang@maastrichtuniversity.nl}\\
\email{** zhoujinhua@ahmu.edu.cn}}
\maketitle              
\begin{abstract}
In this research, we reveal the inborn but hitherto ignored properties of quantitative differential phase contrast (qDPC) imaging: the phase transfer function being an edge detection filter. Inspired by this, we highlight the duality of qDPC between optics and pattern recognition, and propose a simple and effective qDPC reconstruction algorithm, termed Pupil-Driven qDPC (pd-qDPC), to facilitate the phase reconstruction quality for the family of qDPC-based phase reconstruction algorithms. We formed a new cost function in which modified ${L_{0}\text{-norm}}$ was used to represent the pupil-driven edge sparsity, and the qDPC convolution operator is duplicated in the data fidelity term to achieve automatic background removal. Further, we developed the iterative reweighted soft-threshold algorithms based on split Bregman method to solve this modified ${L_{0}\text{-norm}}$ problem. We tested pd-qDPC on both simulated and experimental data and compare against state-of-the-art (SOTA) methods including ${L_{2}\text{-norm}}$, total variation regularization (TV-qDPC), isotropic-qDPC, and Retinex qDPC algorithms. Our model is superior in terms of phase reconstruction quality and implementation efficiency, in which it significantly increases the experimental robustness while maintaining the data fidelity. In general, the pd-qDPC enables high-quality qDPC reconstruction without any modification to the optical system. It simplifies the system complexity and benefits the qDPC community and beyond including but not limited to cell segmentation and PTF learning based on the edge filtering property. 

\end{abstract}
\clearpage
\section{Introduction}\label{section 1}
\begin{figure}
    \centering
    \includegraphics[width = 0.95\textwidth,trim = 0 50 0 50,clip]{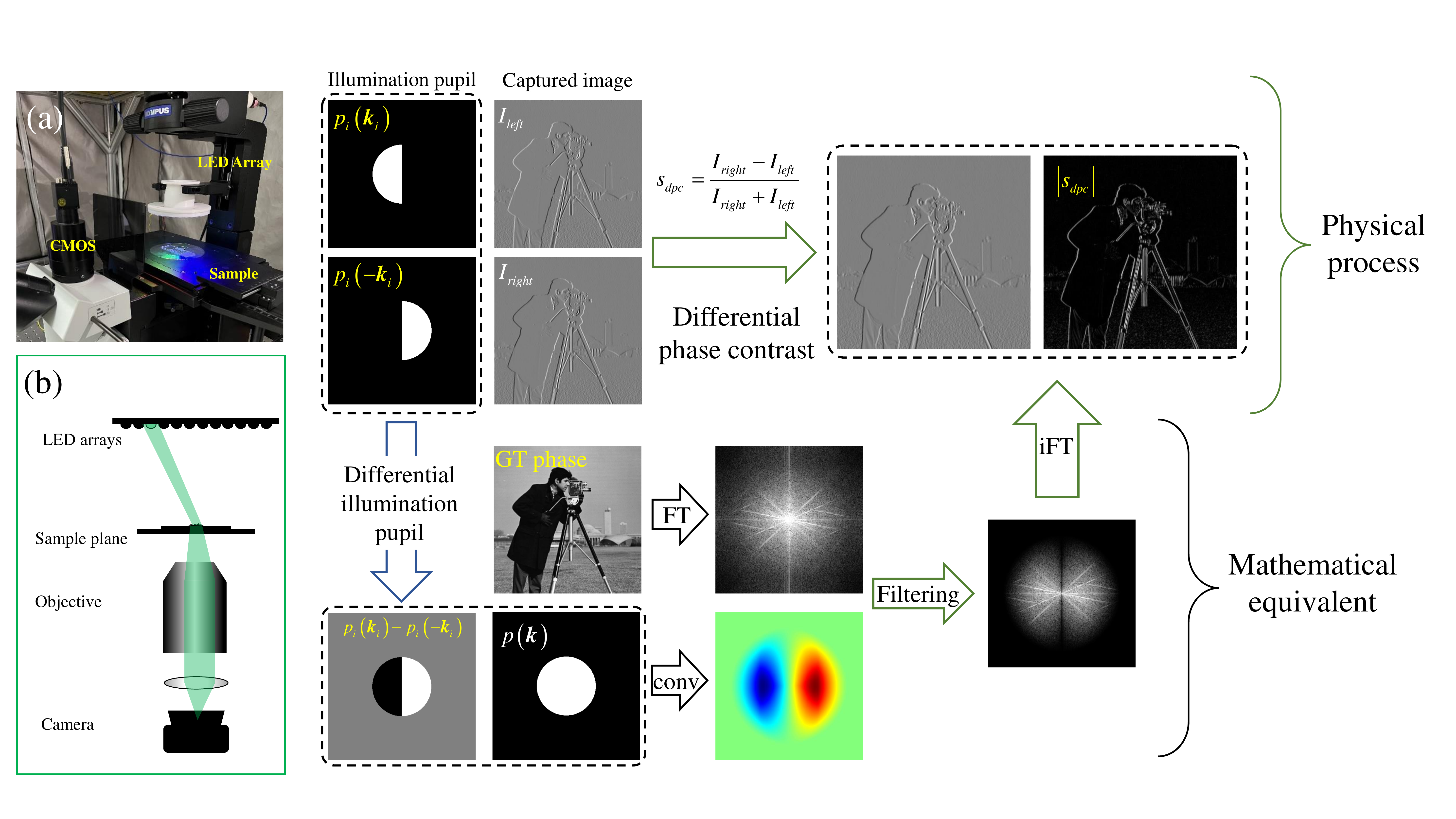}
    \caption{Classical qDPC routine. The sample is illuminated by well-designed oblique illumination pattern to produce phase contrast images. The differential phase contrast (dpc-) image can be obtained by differentiate two phase contrast images under ansymetric oblique illumination. In mathematics, the dpc-image is produced by the Fourier filtering of the Fourier spectrum of sample phase pattern and the phase transfer function, which can be regard as the convolution of differential illumination pupil and the pupil function of the image system}
    \label{Fig.1.}
\end{figure}

Differential phase-contrast microscopy (DPC), a non-interferometric quantitative phase retrieval approach, has been used for label-free and stain-free optical imaging of live biological specimens both in vitro \cite{hamilton1984differential} and in vivo \cite{kandel2019epi,laforest2020transscleral}. A quantitative DPC (qDPC) experimental layout involves a 4-f microscopy system in which a programmable LED or LCD illumination source generates anti-symmetric illumination patterns \cite{tian2015quantitative,lee2015color}. As shown in Fig. \ref{Fig.1.}, with the combination of oblique plane wave illumination and low-pass filtering of the objective lens, the DPC converts the unmeasurable sample phase into a phase-contrast intensity measurements. By collecting at least 4 phase-contrast images with asymmetric illuminations in opposite directions, the phase component of the sample can be reconstructed through a non-blind deconvolution process, where the convolution kernel, in an ideal condition, is defined by the Fourier transform of phase transfer function (PTF) \cite{tian2015quantitative}. The spatial deconvolution is then transformed into Fourier space division.

Despite its many advantages, qDPC is not immune to the imperfections of optical systems since the PTF is ill-conditioned due to  its band-limited \cite{xu2014inverse}, as such even small perturbations in the DPC images can lead to significant errors in the phase reconstruction result, especially under the existence of unavoidable quantization errors and camera noise. The DPC raw images often contain noises and background fluctuations, the latter can result in the appearance of "white-cloud" effects in the phase pattern \cite{zhang2023retinex}, while noise pixels can cause severe ringing artifacts infesting the whole image. Consequently, the applicability of qDPC is limited by its sensitivity to system noise. To improve the performance of qDPC, two families of methods have been proposed: (1) optical system modification, and (2) algorithm modification.

The manuscript is arranged as follows: \textbf{Section \ref{section 2}} briefly reviews off-the-shelf qDPC reconstruction algorithms, while \textbf{Section \ref{section 3}} introduces the insight of our pupil-driven qDPC. In \textbf{Section \ref{section 4}}, we propose the pd-qDPC loss function and a Split Bregman solution. \textbf{Section \ref{section 5}} and \textbf{Section \ref{section 6}} presents experimental results for validation and application of the pd-qDPC, followed by a discussion in \textbf{Section \ref{section 7}} and concluding remarks in \textbf{Section \ref{section 8}}.
\section{Related works}\label{section 2}
\subsection{Improving qDPC by modifying the optical system}
Being an image techniques, the quality of qDPC reconstruction is directly related to the quality of the raw image data. By increasing the quality of data collection, for example, using a high-end camera to increase the signal-to-noise ratio \cite{bonati2022lock}, or using better light source instead of LEDs to suppress the background fluctuation \cite{bonati2020phase}, the final phase deconvolution quality can be improved. However, the noise signal and sample-induced background fluctuation (floating objects in defocusing layers) cannot be completely removed. As such, these methods highly rely on certain experimental conditions and lack of generalization. Also, the high-end hardwares are rather expensive.   

The qDPC reconstruction quality is related to the PTF determined by the illumination pattern and pupil of the objective lens. Since it is inconvenient to directly modify the pupil of the objective lens, the PTF engineering is converted into illumination pattern optimization \cite{chen2018isotropic,lin2018quantitative,fan2019optimal}. In these studies, illumination is optimized using, for example, radially asymmetric patterns \cite{lin2018quantitative}, gradient amplitude patterns \cite{chen2018isotropic}, or ring-shaped sine patterns \cite{fan2019optimal}. As such, the PTF can be better defined and maintain the low frequency components \cite{jenkins2015quantitative}, and the frequency respond of sample phase is improved. However, the problematic background and noise signals are still not solved since they are independent of the optical system. 


\subsection{Improving qDPC by modifying the reconstruction algorithm}
The qDPC phase recovery quality can be also improved by optimization algorithms. The qDPC reconstruction is a maximum a posteriori task in which we would like to find the phase, $\mathbf{\Phi}$, that can best explain the series of collected DPC images \cite{gribonval2011should}. The forward model of qDPC is given by a convolution process that is
\begin{equation}
   \mathbf{S}_n = \mathbf{H}_n \mathbf{\Phi} + \mathbf{\epsilon} \label{(1)}
\end{equation}
where $\mathbf{S}_n$ is the $n^\text{th}$ DPC images ($n$ = 1, 2, ……), which is calculated from the differences between twp phase-contrast images under anti-systematic oblique illuminations
\begin{equation}
   \mathbf{S}_n = \frac{I_{n,l}-I_{n,r}}{I_{n,l}+I_{n,r}}, \label{(2)}
\end{equation}
here, for example, $I_{n,l}$ and $I_{n,r}$ are captured images under anti-systematic oblique illuminations from left and right, respectively. $\mathbf{H}_n$ is the (large) Toeplitz matrix of the point spread function (PSF) denoting the convolution operation between $\mathbf{\Phi}$ and the PSF. $\mathbf{\epsilon}$ denotes the noise signal. The phase reconstruction for $\mathbf{\Phi}$ is given by an optimization problem where
\begin{equation}
   \mathbf{\Phi} = \text{argmin}_{\mathbf{\Phi}} F \left(\mathbf{S}_n,\mathbf{\Phi}\right)+\alpha R(\mathbf{\Phi}), \label{(3)}
\end{equation}
where $F \left(\mathbf{S}_n; \mathbf{\Phi}\right)$ denotes the likelihood, also known as data fidelity term that measures the distance between observation $\mathbf{S}_n$ and mathematical prediction $\mathbf{H}_n \mathbf{\Phi}$. $R(\mathbf{\Phi})$ is the regularization/penalty term that restrict the solution of $\mathbf{\Phi}$. $\alpha$ is the penalty parameter that controls the strength of the regularization. Functions $F$ and $R$ are summarized in Tab. \ref{Table 1} according to published literature.

\begin{table}
\renewcommand\arraystretch{2}
\caption{List of qDPC solver. $\alpha$ and $\beta$ are penalty parameters. $ \mathcal{G} $ denotes convolution with Gaussian kernel.}
\label{Table 1}
\begin{adjustbox}{width=\columnwidth,center}
\begin{tabular}{c r c l l}
\hline
 &Data fidelity term &+ &Penalty term &Prior\\
\hline
\multirow{3}{*}{Traditional}   &\multirow{3}{*}{$\sum_{n=1}^{N}\left \|\mathbf{S}_n - \mathbf{H}_n\mathbf{\Phi } \right \|_2^2 $} &\multirow{6}{*}{+}&$\alpha\left\|\mathbf{\Phi}\right\|_2^{2}$ &$L_2\text{-norm} $ \cite{tian2015quantitative}\\
& & &$\alpha\left\|\nabla\mathbf{\Phi}\right\|_1$ &Total variation (TV) \cite{rudin1992nonlinear}\\
& & &$\alpha\left\|\nabla\nabla\mathbf{\Phi}\right\|_2^{2}$ &High-order TV \cite{chan2000high}\\
\multirow{3}{*}{Retinex-qDPC \cite{zhang2023retinex}}&\multirow{3}{*}{$\sum_{n=1}^{N}\left \|\nabla\mathbf{S}_n - \mathbf{H}_n\nabla\mathbf{\Phi } \right \|_2^2 $}& &$\alpha\left\|\mathbf{\Phi}\right\|_2^{2} + \beta\left\|\mathcal{G}\mathbf{\Phi}\right\|_2^{2}$ & Isotropic qDPC \cite{chen2018isotropic}\\
& & &$\alpha\left\|\mathbf{\Phi}\right\|_0$ &Sparse-sample qDPC \cite{peng2023quantitative}\\
& & & Other penalty &\\
\hline
\end{tabular}
\end{adjustbox}
\end{table}

Traditional qDPC uses the $L_2$-norm to measure the distance between $\mathbf{S}_n$ and $\mathbf{H}_n \mathbf{\Phi}$, assuming that both noise is sampled from a Gaussian distribution. However, this approach is not effective at suppressing background fluctuations. To address this, the Retinex-qDPC \cite{zhang2023retinex} was proposed, which uses the image spatial gradient for data fidelity. Since the gradient of image is sparse, the $L_1\text{-norm}$ is further adapted to promote the sparsity forming the $L_1\text{-retinex}$-qDPC.  

For the penalty term, the $L_2$-norm regularization is most commonly used in DPC experiments, which implies that the phase is also sampled from Gaussian distribution \cite{rennie2003l2}. However, the penalty parameter $\alpha$  should be carefully selected as it may suppress the value of entire $\mathbf{\Phi}$ \cite{zhang2018three} and introduce black-halo effect to the reconstructed phase pattern \cite{zosso2013unifying}. The total variation \cite{rudin1992nonlinear} and its high-order version \cite{chan2000high} suppress the noise signal by eliminating the small value of image gradient (higher-order gradient). For qDPC experiment, the TV-regularization can significantly suppress the noise signal \cite{chen2018quantitative}, but using first-order TV regularization may lead to loss of small structures and a stair-step effect.

In \cite{chen2018isotropic}, the gradient Tikhonov regularization was employed and a Gaussian filter is used to estimate the background. By minimizing the $L_2\text{-norm}$ of spatial gradient of the phase pattern and the Gaussian filtered pattern, background fluctuation and noise signals can be suppressed. However, mild white-cloud effect and ringing effect may still be present \cite{chen2018isotropic}. Recent publication on qDPC using the deep-image-prior \cite{Chen:23} shows promising results on denoising for qDPC reconstruction, while still cannot handle the intensity mismatch problem. It also share the same shortage of deep-image-prior \cite{ulyanov2018deep} where the high-frequency structures is hard to be recovered.  

Although the above-mentioned methods are widely used for solving qDPC inverse problem, they are not fully optimized according to specific feature of qDPC raw data, and the qDPC community are seeking for a methods that can simultaneously address the background fluctuation and noise problem for user-friendly and stable application of qDPC. 

\subsection{Our contribution}
In this research, we proposed the Pupil-driven qDPC (pd-qDPC) for the phase reconstruction. The pd-qDPC modifies the cost function of qDPC deconvolution using the PTF information determined by the optical system, to significantly improve the reconstruction quality, and address both noise signals and background fluctuations. This algorithm is suitable for all DPC experiments and does not require any modifications to the optical devices. Figure 2 demonstrates the effectiveness of the pd-qDPC phase reconstruction results. 

Our contributions are summarized as follows: 
1)	we analysis the spectrum and spatial response of PTF in qDPC experiment and found that the \textbf{PTF works identically as an edge detection filter}. With that, we are able to add the edge sparsity prior as a penalty term for noise suppressing.

2)	Being an edge detection filter, we are able to duplicate the DPC convolution in the data fidelity term to achieve a Retinex-like filtering. \textbf{This allows the automatically background rectified qDPC reconstruction.}

3)	With (1) and (2), \textbf{we propose a novel cost function termed pd-qDPC cost for quantitative phase reconstruction task.} The pd-qDPC cost is solved using split Bregman methods with our developed iterative reweighted soft-threshold operator.

4)	\textbf{Our method performs well on both synthetic DPC data and real DPC experiments, and outperforms state-of-the-art (SOTA) algorithms.} 

5)	The illumination pupil of qDPC can be designed and learning from the edge response cost function of given edge pattern. 

In the following section we introduce the concept of pupil-driven qDPC (pd-qDPC) and prove why it works mathematically. 
\section{Phase transfer function as edge filters}\label{section 3}
The pd-qDPC is based on a proposition of the PSF of qDPC experiment: the PSF, generated by the Fourier transform of PTF, is a band-limited edge detection filter (twice the size of $NA/\lambda$, NA is the numerical aperture of the image system, $\lambda$ is the wavelength). 

Why this is hold for qDPC? Intuitivly, the phase difference of two oblique illuminated images happens on the area where the phase pattern changes rapidly. The differential operation of two phase contrast image in ideal condition should completely cancel out the background components since they are precisely matched pixel-wisely. Edge structure of the phase pattern can therefore be observed. Similar phenomenon can be also observed in the differential operation of transport of intensity \cite{beleggia2004transport}.

\begin{figure}
    \centering
    \includegraphics[width = 0.98\textwidth,trim = 0 20 590 0,clip]{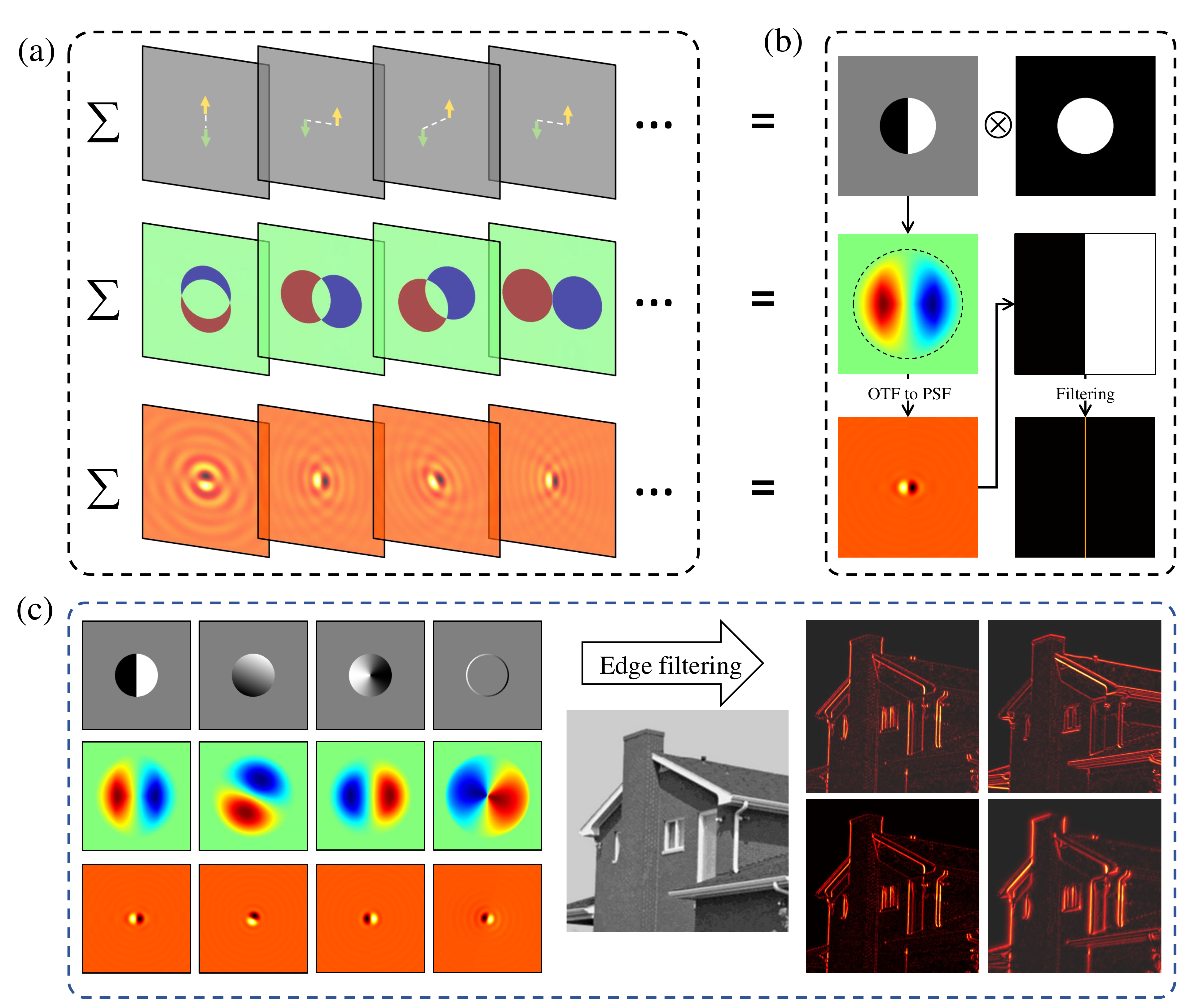}
    \caption{Decomposition of illumination pupil function $\mathbf{Q}_n$. (a) In first row, $\mathbf{Q}_n$ is expressed as superposition of dual-pulse function of different distance and direction. The yellow arrow denotes positive pulse, and the green arrow denotes the negative pulse. The second row denotes the convolution between $\mathbf{Q}_n$ and $\mathbf{P}$. The third row denotes the corresponding PSF. (b) and (c) show the edge detection feature of psf/PTF for different PTF layouts, and images. The psf is zoomed-in for visual performance.}
    \label{Fig.2.}
\end{figure}

To mathematically illustrate this observation, we consider the property of point spread function $\mathbf{H}_n$ of the differential phase contrast image $\mathbf{S}_n$ in Fourier domain, which is given by 
\begin{equation}
   \mathbf{H}_n(\mathbf{r})=\mathbf{F}^{-1}\widetilde{\mathbf{H}}_n(\mathbf{f}), \label{(4)}
\end{equation}
where $\mathbf{r}=(x,y)$ is the spatial coordinates. $\mathbf{F}^{-1}$ is the inverse Fourier transform. $\mathbf{f}=\left(f_x,f_y\right)$ is the coordinate vector in Fourier space. $\widetilde{\mathbf{H}}_n$ is the  phase transform function (PTF). For qDPC in an ideal condition, $\widetilde{\mathbf{H}}_n$ is determined by \cite{tian2015quantitative}:
\begin{equation}
   \widetilde{H}_n(\mathbf{f})=iA\iint q_n(\bm{\uprho})P^*(\mathbf{\bm{\uprho}})
\left [ P(\bm{\uprho}+\mathbf{f}) - P(\bm{\uprho}-\mathbf{f})\right ] \text{d}^2\bm{\uprho}, \label{(5)}
\end{equation}
where $A$ is an energy normalization coefficient. $q_n$ is the intensity of the incident plane wave. $P$ is the aberration-free pupil function of the objective lens which is given by
\begin{equation}
P(\mathbf{f})=\text{circ}\left(\frac{\lambda}{NA}\left | \mathbf{f} \right | \right)=
\begin{cases}
 1 & \text{ if } \left | \mathbf{f} \right |<NA/\lambda \\
 0 & \text{ if } \left | \mathbf{f} \right |=elsewise
\end{cases}
,\ \left | \mathbf{f} \right |=\sqrt{f_x^2+f_y^2}. \label{(6)}
\end{equation}
we rewrite the third term in Eq. (\ref{(5)}) as the convolution between $P$ and $\delta$ function with is 
\begin{equation}
\begin{aligned}
P(\bm{\uprho}+\mathbf{f})-P(\bm{\uprho}-\mathbf{f}) & = P(\mathbf{f})\otimes \left [ \delta(\mathbf{f}+\bm{\uprho})-\delta(\mathbf{f}-\bm{\uprho})  \right ] \\
& = \iint P(\bm{\uptau})\left [ \delta(\mathbf{f}+\bm{\uprho}-\bm{\uptau})-\delta(\mathbf{f}-\bm{\uprho}-\bm{\uptau})  \right ]\text{d}^2\bm{\uptau}
\end{aligned}, \label{(7)}
\end{equation}
and $\otimes$ denotes the 2D convolution. Submitting Eq. (\ref{(7)}) into Eq. (\ref{(5)}) and perform integration w.r.t. $\bm{\uprho}$ first we obtain
\begin{equation}
\begin{aligned}
\widetilde{H}_n\left ( \mathbf{f} \right ) & = iA\ P(\mathbf{f}) \otimes Q_n(\mathbf{f}) \\
&\text{where}\ \  Q_n(\mathbf{f})=\int_{0}^{\frac{NA}{\lambda} }\int_{\theta_0-\frac{\pi}{2}}^{\theta_0+\frac{\pi}{2}}  q_n(\bm{\uprho})P^*(\bm{\uprho})\left [ \delta(\mathbf{f}+\bm{\uprho})-\delta(\mathbf{f}-\bm{\uprho})  \right ] \rho \text{d}\rho\text{d}\theta
\end{aligned}. \label{(8)}
\end{equation}
The integral area in Eq. (\ref{(8)}) is determined according to band-limited feature of $P$. $Q_n$ in Eq. (\ref{(8)}) is an odd function that denotes the pattern of illumination pupil being the superposition two groups of plane waves. In Fourier domain, each plane wave is given by a delta function as shown in Fig. \ref{Fig.2.} (a), propagating in direction of $(\rho \cos{\theta},\rho \sin{\theta},\sqrt{1-\rho^2})$. 

$\theta_0$ controls the central illumination direction of the two groups of plane waves, and for example, when $q_n (\mathbf{\bm{\uprho}}) = 1$ and $\theta_0=0$, Eq. (\ref{(8)}) generates two half-circles one is positive on the right-hand side of y-axis, and the other one is negative on the left hand side of y-axis as shown in Fig. \ref{Fig.2.} (b). Since convolution and the integration can be switched, the PTF can be regarded as the superposition of convolution results of dual-delta function and $P$ which is also listed in Fig. \ref{Fig.2.} (a). With Eq. (\ref{(5)}) to Eq. (\ref{(8)}), we have the following 


\begin{lemma} \label{Lemma1}
When $\theta_0 = 0$ the PSF described by Eq. (\ref{(8)}) is an edge detection filter for edge drawing along y-axis, since it produces the highest responce along y-axis when it is filtering with a unit-step function $u(x,y) = \text{sign} (x)$. 
\end{lemma}

\begin{proof}
One has the Fourier spectrum of unit-step function:
\begin{equation}
\widetilde{u} (f_x,f_y)=\begin{cases}
\frac{1}{i f_x}   & f_y = 0\\
 0 & elsewise
\end{cases}. \nonumber
\end{equation}
Combined with (\ref{(5)}) to Eq. (\ref{(7)}) and performing inverse Fourier transform yielding the expression of filter responce $\left|g(x,y)\right|$
\begin{equation}
    \left|g(x,y)\right| = \left| H_n \otimes u \right| = \left| A\int_{0}^{\frac{NA}{\lambda} } \left|\frac{\left \{ P(\mathbf{f} )\otimes Q_n(\mathbf{f} ) \right \} |_{f_y=0}}{f_x} \right|
\cos \left ( 2 \pi f_x x \right ) \text{d} f_x  \right| \nonumber
\end{equation} 
note that $P(\mathbf{f} )\otimes Q_n(\mathbf{f})  < 0 $ for $f_x \in [0,NA/\lambda]$. Following triangle inequality we have
\begin{equation} 
\left|g(x,y)\right| < A \int_{0}^{\frac{NA}{\lambda} } \left|\frac{\left \{ P(\mathbf{f} )\otimes Q_n(\mathbf{f} ) \right \} |_{f_y=0}}{f_x} \right|
\left|\cos \left ( 2 \pi f_x x \right )\right| \text{d} f_x \le \left|g(x \equiv 0,y)\right|\nonumber
\end{equation} 
As such, the filter-responce strength becomes the strongest at $x=0$ (along y-axis).

\end{proof}

With Lemma \ref{Lemma1}, we have the following 

\begin{lemma} \label{Lemma2}
For an arbitray $\theta_0=\theta$, the PSF described by Eq. (\ref{(8)}) is an edge detection filter for edge drawing along vector $\mathbf{v}=(-\sin{\theta_c},\cos{\theta_c})$.
\end{lemma}
\begin{proof}
This situation can be deduced to Lemma \ref{Lemma1} if one rotates the PSF as well as the edge around the origin for $\theta$ angle. Then, following Lemma \ref{Lemma1} this is proved. The PSF's edge detection is shown in Fig. \ref{Fig.2.} (b) and (c) for different input images.
\end{proof}


Given the edge-detection feature of qDPC's PSF, the PTF itself can be used to regularize the solution of qDPC inverse problem, and boosting the phase reconstruction results. As such, we propose the pupil-driven qDPC (pd-qDPC) reconstruction routine to achieve high experimental robustness quantitative phase recovery. 
\section{Pupil-driven qDPC}\label{section 4}
\subsection{pd-qDPC inverse problem}
In this section, we put forward a concrete qDPC deconvolution model and one effective optimization scheme. In the framework of Maximum likelihood estimation (MAP) in Eq. (\ref{(3)}), the regularization term $R(\mathbf{\Phi})$ implies the prior distribution of parameter $\mathbf{\Phi}$. Since image edge follows Heavy-tailed distribution \cite{10.1007/978-3-642-40246-3_8} and $\left |\mathbf{H}_n \mathbf{\Phi} \right |$ is a group of edge extracted from $\mathbf{\Phi}$, we are able to adapt the Heavy-tailed distribution to $\left |\mathbf{H}_n \mathbf{\Phi} \right |$ to form our pupil-driven edge regularization term. Further, we use $\mathcal{P} (x)=\mathrm{exp} \left [ \mathrm{exp} \left ( \left | \omega x \right | \right ) -1  \right ]$ to approximate the Heavy-tailed distribution, and $\omega$ controls the width of the curve. In addition, the first order spatial gradient of $R(\mathbf{\Phi})$ also follows the Heavy-tailed distribution, and we use the same approximation to treat this. In general, the regularization term of our pupil-driven qDPC is given as
\begin{equation}
R\left ( \mathbf{\Phi}   \right ) =\alpha \sum_{n=1}^{N}\left \| 1-\mathrm{exp}\left ( -\omega  \left | \mathbf{H}_n \mathbf{\Phi} \right |  \right )   \right \|_1
+ \beta \left \| 1- \mathrm{exp}\left (-\omega \left |\nabla \mathbf{\Phi} \right| \right ) \right \|_1. \label{(9)}
\end{equation}
The first term in Eq. (\ref{(9)}) denotes our proposed pupil-driven edge regularization term, and the second term denotes the gradient-edge regularization term. 

Further, based on the edge-detection feature of PTF, we developed the pupil-driven fidelity term by applying $\mathbf{H}_n$ to both $\mathbf{S}_n$ and $\mathbf{H}_n\mathbf{\Phi}$ again, which is 
\begin{equation}
F=\sum_{n=1}^{N} \left \| \mathbf{H}_n\mathbf{S}_n - \mathbf{H}_n \mathbf{H}_n \mathbf{\Phi}    \right \|_2^2. \label{(10)}
\end{equation}
Since $\mathbf{H}_n$ is an edge filter, Eq. (\ref{(10)}) shares the similar idea with Retinex-qDPC in which the background fluctuation is suppressed by using the edge structures \cite{zhang2023retinex}. Although $\mathbf{H}_n$ is band-limited, using them to filter the image again will not cause lose of information since the spectral information of $\mathbf{S}_n$ in an ideal condition was already truncated by $\mathbf{H}_n$ during the optical process of data collection. Further, $\mathbf{H}_n$ also helps erasing low-frequency components (maybe) caused by background mismatch that should not appear to the ideal DPC images.

Combining Eq. (\ref{(9)}) and Eq. (\ref{(10)}), the pupil-driven qDPC cost function is given as 
\begin{equation}
    \mathcal{L} \left ( \mathbf{\Phi} \right ) =\sum_{n=1}^{N} \left \| \mathbf{H}_n\mathbf{S}_n - \mathbf{H}^2_n \mathbf{\Phi}    \right \|_2^2
+\alpha \sum_{n=1}^{N}\left \| f\left ( \mathbf{H}_n \mathbf{\Phi}  \right )   \right \|_1
+ \beta \left \| f\left (  \nabla \mathbf{\Phi}  \right ) \right \|_1, \label{(11)}
\end{equation}
where $f(x)=1-\text{exp}(- \omega \left| x \right|)$. The solution of $\mathbf{\Phi}$ is given by minimize the cost function, $\mathbf{\Phi} = \text{argmin} \   \mathcal{L} \left ( \mathbf{\Phi} \right )$.
\subsection{Split-Bregman solver}
We use the split Bregman method \cite{goldstein2009split} to solve the two non-convex $L_1 \text{-norm}$ involved problems in Eq. (\ref{(11)}). By introducing $N + 2$ auxiliary variables $\mathbf{\Psi}_n, (n=1,2, \cdots,  N)$ and $\mathbf{G}$ with respect to $\mathbf{H}_n \mathbf{\Phi}$ and $\nabla \mathbf{\Phi}$, Eq. (\ref{(11)}) is converted into $N+3$ subproblems. At t-th iteration, these sub-problems are

\begin{equation}
\left\{\begin{matrix}
\begin{aligned}
    \mathbf{\Phi }^t & = \mathrm{argmin} \sum_{n = 1}^{N} \left \| \mathbf{H}_n\mathbf{S}_n - \mathbf{H}^2_n \mathbf{\Phi}    \right \|_2^2
+\alpha_0\sum_{n  = 1}^{N} \left \| \mathbf{H}_n \mathbf{\Phi} -\mathbf{\Psi }_n^{t-1}-\mathbf{b}_n^{t-1}    \right \|_2^2 \\
  & \ \ \ \ \ \ \ \ \ \ \ \ \ \ \ \ \ \ \ \ \ \ \ \ \ \ \ \ \ \ \ \ \ \ \ \ \ \ \ \ \ \ \ \ \ \ \ \ \ \ \ \ \ \ +  \beta_0 \left \| \nabla \mathbf{\Phi} - \mathbf{G}^{t-1} - \mathbf{d} ^{t-1} \right \|_2^2  \\
\mathbf{\Psi }_n^{t} & = \mathrm{argmin}\left \| \mathbf{H}_n \mathbf{\Phi}^t -\mathbf{\Psi }_n^{t-1}-\mathbf{b}_n^{t-1}    \right \|_2^2
+\frac{\alpha}{\alpha_0}\left \| f\left (  \mathbf{\Psi }_n\right )  \right \|_1,\ n = 1,\ 2,\ \cdots,\  N  \\
\mathbf{G }^{t} & = \mathrm{argmin}\left \| \nabla \mathbf{\Phi} - \mathbf{G}^{t-1} - \mathbf{d} ^{t-1} \right \|_2^2
+\frac{\beta}{\beta_0}\left \| f\left (  \mathbf{G} \right )  \right \|_1  \\
\mathbf{b}_n^{t} & = \mathbf{b}_n^{t-1}+\mathbf{\Psi }_n^t-\mathbf{H}_n \mathbf{\Phi}^t \\
\mathbf{d}^{t} & = \mathbf{d}^{t-1}+\mathbf{G }^t-\nabla \mathbf{\Phi} ^t
\end{aligned}
\end{matrix}\right. \label{(12)}
\end{equation}
where $\mathbf{\Psi }_n^0 = \mathbf{b}_n^0 = 0$, $\mathbf{b}_n^t$ and $\mathbf{d}^t$ is the Bregman parameters. $\beta_0$ is the penalty parameter for the quadratic term and $\alpha_0=\beta_0=1$ throughout the manuscript. 

The $\mathbf{\Phi }$-subproblem in Eq. (\ref{(12)}) is pure-quadratic and can be solved by setting the derivative with respect to $\mathbf{\Phi }$ to zero. The closed-form solution is given by
\begin{equation}
    \mathbf{\Phi }=\frac{\sum_{n=1}^{N} \left [ \mathbf{H}_n^\mathsf{T}\mathbf{H}_n^\mathsf{T}\mathbf{H}_n \mathbf{S}_n +\alpha _0\mathbf{H}_n^\mathsf{T}\left ( \mathbf{\Psi}_n^{t-1}  +\mathbf{b}_n^{t-1} \right ) \right ]  + \beta_0  \nabla ^\mathsf{T} \left( \mathbf{G}^{t-1} + \mathbf{d}^{t-1} \right)  }{
 \sum_{n=1}^{N} \left (\mathbf{H}_n^\mathsf{T}\mathbf{H}_n^\mathsf{T}\mathbf{H}_n\mathbf{H}_n + \alpha _0\mathbf{H}_n^\mathsf{T}\mathbf{H}_n \right ) + \beta_0 \nabla^\mathsf{T}\nabla  +\eta } \label{(13)}
\end{equation}
where $\mathsf{T}$ denotes the transpose of matrix. Eq. (\ref{(13)}) can be accelerated using Fourier transform. 

The $\mathbf{\Psi }$ sub-problem can be solved using reweighted soft-threshold methods \cite{ochs2015iteratively}, where 
\begin{equation}
    \mathbf{\Psi } _n^t=\text{sign}\left ( \mathbf{H}_n\mathbf{\Phi } -\mathbf{b}_n^{t-1} \right )
 \max \left \{ \left | \mathbf{H}_n\mathbf{\Phi } -\mathbf{b}_n^{t-1} \right | - \frac{\alpha}{\alpha_0}
\frac{\left | {f}' \left ( \mathbf{H}_n\mathbf{\Phi } -\mathbf{b}_n^{t-1} \right )  \right | }{\left | {f}' \left ( \alpha/\alpha_0  \right )  \right | }  ,\ 0  \right \}, \label{(14)}
\end{equation}
for anisotropic manner, and 
\begin{equation}
    \mathbf{\Psi } _n^t=\frac{\mathbf{H}_n\mathbf{\Phi } -\mathbf{b}_n^{t-1}}{\mathbf{\Lambda}}  
 \max \left \{ \mathbf{\Lambda} - \frac{\alpha}{\alpha_0}
\frac{\left | {f}' \left ( \mathbf{\Lambda} \right )  \right | }
{\left | {f}' \left ( \alpha/\alpha_0  \right )  \right | }  ,\ 0  \right \},
\ \mathbf{\Lambda}=\sqrt{\sum_{n=1}^{N}\left ( \mathbf{H} _n \mathbf{\Phi}-\mathbf{b}_n^{t-1}\right )^2  } . \label{(15)}
\end{equation}
for isotropic manner. 

\begin{figure}
    \centering
    \includegraphics[width = 0.75\textwidth,trim = 0 150 0 150,clip]{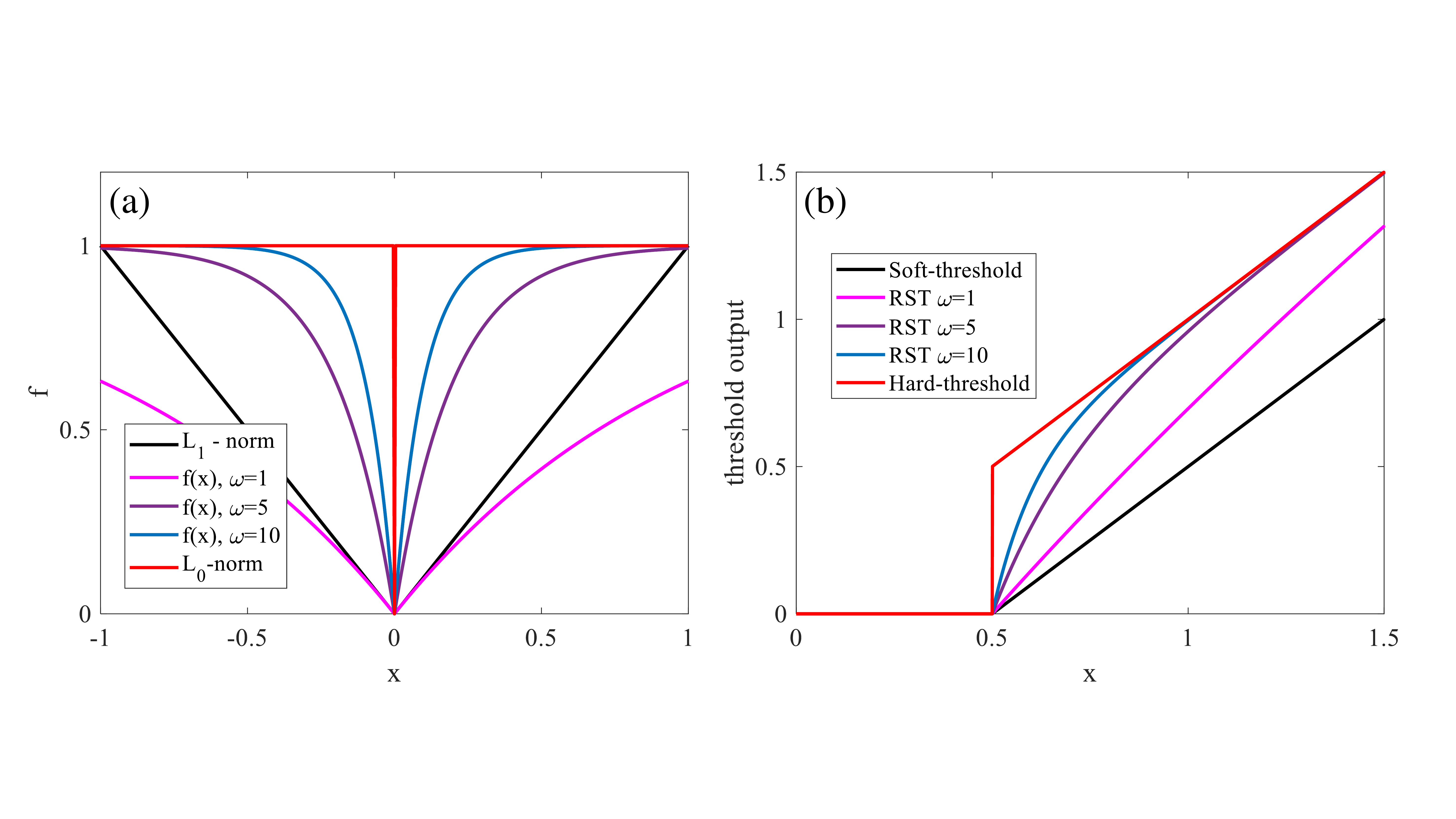}
    \caption{Demonstration of function $f$ with respect to different value of $\omega$ and its RST output. (a) curves for function $f$, $L_1\text{-norm}$, and $L_0\text{-norm}$ (b) threshold operation output. The soft and hard threshold operations are solution for $L_1\text{-norm}$ and $L_0\text{-norm}$ problem}
    \label{Fig.4.}
\end{figure}

Similar procedure applied to solve the $\mathbf{G}$ sub-problem by replacing symbols $\mathbf{\Psi }$ and $\mathbf{H}_n$ in Eqs. (\ref{(14)}), (\ref{(15)}) with $\mathbf{G}$ and $\nabla$, respectively. The reweight soft-threshold in Eq. (\ref{(14)}) and Eq. (\ref{(15)}) can behave between the hard-threshold operation and soft-threshold operation by adjusting parameter $\omega$ as shown in Fig. \ref{Fig.4.}. In general, the algorithm for pd-qDPC is listed in \ref{Table Algorithm 1}.

\begin{table}
\renewcommand\arraystretch{1.2}
\caption{Algorithm 1, Pupil-driven qDPC phase deconvolution}
\label{Table Algorithm 1}
\begin{adjustbox}{width=0.6\columnwidth,center}
\begin{tabular}{l}
\hline
\textbf{Input:} Raw PDC-images $\mathbf{S}_n$, $\mathbf{H}_n$, Penalty parameters $\alpha$, $\beta$      \\ \hline
Initializing   $\mathbf{\Psi}_n^0 = \mathbf{G}^0 = \mathbf{b}_n^0 = \mathbf{d}_n^0 = 0$                                  \\
\textbf{While} $t < loop_{max}$ \textbf{do}               \\
\ \ \ \ \ Solving $\mathbf{\Phi}$ using Eq. (13)             \\
\ \ \ \ \ Solving $\mathbf{\Psi}_n^t$ using Eq. (14) or (15)     \\
\ \ \ \ \ Solving $\mathbf{G}^t$ similar to $\mathbf{\Psi}$                   \\
\ \ \ \ \ Updating $\mathbf{b}_n^t$ and $\mathbf{d}^t$ using Eq. (12)                    \\
\ \ \ \ \ $t = t + 1$                                          \\
\textbf{End}                                                \\
\textbf{Output}: Phase $\mathbf{\Phi}$                      \\ \hline
\end{tabular}
\end{adjustbox}
\end{table}

\subsection{Adaptive noise sensor}
The pd-qDPC model Eq. (\ref{(10)}) has two penalty parameters, $\alpha$ and $\beta$, that need to be manually adjusted according to the noise level of the input images $\mathbf{S}_n$ which is band-limited ideally. Inspired by this, we designed an adaptive noise sensor that evaluate the noise level from the Fourier domain that locates beyond the band-limited ($ |\mathbf{f}|> 2NA/\lambda$). Given the Laplacian operator $\mathbf{\Gamma} =\left ( -1,2,-1;2,-4,2;-1,2,-1  \right ) $, the noise-level of input $\mathbf{S}$ is depicted as \cite{immerkaer1996fast}
\begin{equation}
    \sigma =\frac{1}{5}\sqrt{\frac{\pi}{2}} \frac{1}{NWH}  
\sum_{n=1}^{N}\sum_{x=1}^{W}\sum_{y=1}^{H}  \left | \text{High-pass-filter}_{|\mathbf{f}| > \frac{2NA}{\lambda} } \left ( \mathbf{S}_n \right )    \otimes \mathbf{\Gamma}     \right |. \label{(16)}
\end{equation}
The high-pass filter denotes filtering the image and erasing all frequency components that within the band-limit. For noise corrupted images, Laplacian operator can enlarge the noise signals as it is sensitive to noise pixels, and the average absolute pixel value of the filtered image can be regarded as a measurement of noise level. 

Accordingly, the parameters $\alpha$ and $\beta$ can be attentively given by
\begin{equation}
    \left\{\begin{matrix}
\begin{aligned}
\alpha&=\sigma/2  \\
\beta&=\sigma/10
\end{aligned}
\end{matrix}\right.\label{(17)}
\end{equation}
\section{Experimental results}\label{section 5}
\subsection{Simulation study}

\begin{figure}
    \centering
    \includegraphics[width = 0.95\textwidth,trim = 0 1160 0 0,clip]{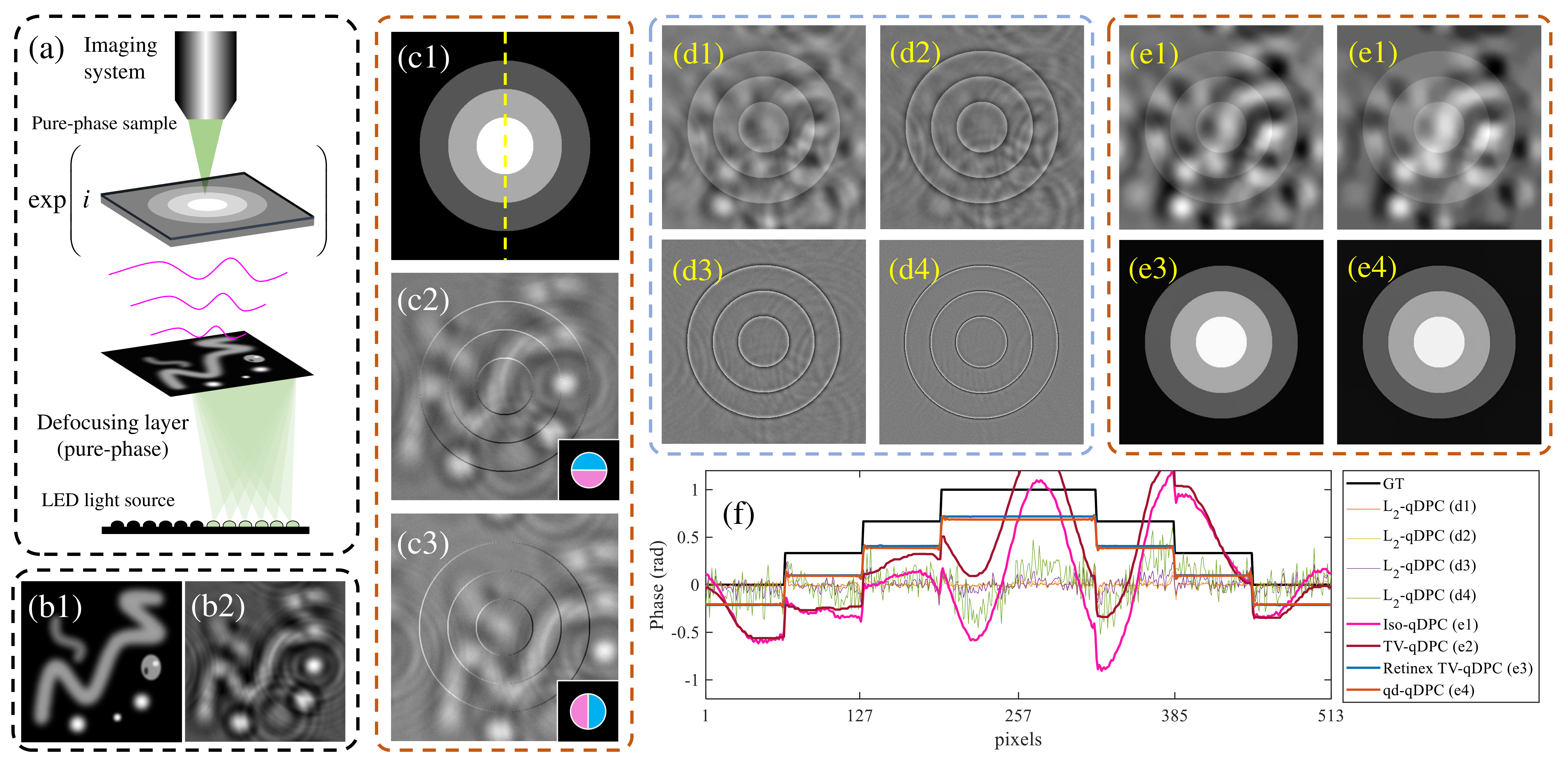}
    \caption{Simulation study on different qDPC-reconstruction methods with defocusing layers. (a) Sketch of optical layout. (b1) is the simulated pure-phase objects in deeper floating layers, $z = -1$ cm, and (b2) is the intensity pattern when it is propagated (diffracted) to the current focal plane ($z = 0$). (c1) is the phase ground truth (GT) of the pure-phase sample.  (c2) and (c3) are differential phase contrast images corrupted by the shadows of the defocal layer. (d1) to (d4) are reconstructed phase using $L_2\text{-qDPC}$ with $\alpha=0.0001$, $\alpha=0.001$, $\alpha=0.01$, $\alpha=0.1$, respectively. (e1) to (e4) are results of Iso-qDPC, TV-qDPC, Retinex TV-qDPC, and our proposed pd-qDPC. (f) is the quantitative phase profile.}
    \label{Fig.5.}
\end{figure}

First, we conducted simulation study on simulated qDPC imaging when floating defocus layer appears in Fig. \ref{Fig.5.} (a). The phase pattern of the pure-phase object (wedding cake) is shown in Fig. \ref{Fig.5.} (c1). The brightest disk denotes 1 rad, while the outer rings denotes 2/3 and 1/3 rads. The groundtruth (also pure-phase) of the defocusing layer is shown in Fig. \ref{Fig.5.} (b1), and is propagated to the current on-focuse layer and the intensity, as shown in Fig. \ref{Fig.5.} (b2), to mimics the background fluctuation caused by defocusing pattern. The NA of this system is 0.25 and so does illumination NA. The wavelength $\lambda=0.532 \ \upmu\text{m}$. The system magnification is $\times 10$, and the pixel size of the camera is $4 \upmu\text{m}$. We also added strong Gaussian noise to the DPC images so that the signal-to-noise (SNR) is 5 db.

The pd-qDPC was compared against other 4 methods namely: $L_2\text{-qDPC}$, Iso-qDPC, TV-qDPC, and Retinex TV-qDPC (TV-qDPC with Retinex data fidelity term). Within, the $L_2\text{-qDPC}$ was tested with variate values of regularization strength as shown in Fig. \ref{Fig.5.} (c1) and Fig. \ref{Fig.5.} (c4). The value of other penalty parameters was determined using Eq. (\ref{(17)}) for remain methods. According to Fig. \ref{Fig.5.}, the $L_2\text{-qDPC}$ fails to recover the phase information due to large noise signals and presents of background fluctuation, and so does Iso-qDPC and TV-qDPC. The phase information is severely degraded due to the shadows. On the contrary, both Retinex TV-qDPC and pd-qDPC show their ability on background rectifying. Further as shown in Fig. \ref{Fig.5.} (e), the recovered phase is consistent with the ground truth (GT) with only constant shift. 

\begin{table}
\renewcommand\arraystretch{1.5}
\caption{Comparison study on simulated data with background mismatch and different strength of noise.}
\label{Table 2}
\begin{adjustbox}{width=\columnwidth,center}
\begin{tabular}{cccccccccc}
\hline
\multirow{3}{*}{\textbf{Metrics}} & \multirow{3}{*}{SNR} & \multicolumn{8}{c}{Methods}                                                                                                                                                                  \\ \cline{3-10} 
                                  &                      & \multicolumn{4}{c}{$L_2$-qDPC}                                          & \multirow{2}{*}{Iso-qDPC} & \multirow{2}{*}{TV-qDPC} & \multirow{2}{*}{Retinex TV-qDPC} & \multirow{2}{*}{pd-qDPC} \\ \cline{3-6}
                                  &                      & $\alpha = 0.1$ & $\alpha = 0.01$ & $\alpha = 0.001$ & $\alpha = 0.0001$ &                           &                          &                                  &                          \\ \hline
\multirow{3}{*}{rpSNR (↑)}        & 5                    & 2.027 ± 0.000  & 2.036 ± 0.000   & 1.947 ± 0.001    & 1.053 ± 0.008     & -6.974 ± 0.003            & -13.230 ± 0.004          & 11.745 ± 0.088                   & \textbf{20.325 ± 0.035}  \\
                                  & 10                   & 2.031 ± 0.000  & 2.057 ± 0.000   & 2.045 ± 0.001    & 1.445 ± 0.003     & -6.974 ± 0.002            & -13.223 ± 0.002          & 15.954 ± 0.092                   & \textbf{20.735 ± 0.049}  \\
                                  & 15                   & 2.032 ± 0.000  & 2.063 ± 0.000   & 2.077 ± 0.000    & 1.576 ± 0.001     & -6.974 ± 0.001            & -13.223 ± 0.001          & 18.845 ± 0.044                   & \textbf{20.814 ± 0.028}  \\
\multirow{3}{*}{PSNR (↑)}         & 5                    & 8.523 ± 0.000  & 8.529 ± 0.000   & 8.472 ± 0.001    & 7.884 ± 0.006     & 1.219 ± 0.003             & -4.802 ± 0.004           & 12.114 ± 0.014                   & \textbf{12.731 ± 0.001}  \\
                                  & 10                   & 8.525 ± 0.000  & 8.542 ± 0.000   & 8.534 ± 0.000    & 8.147 ± 0.002     & 1.219 ± 0.001             & -4.796 ± 0.002           & 12.550 ± 0.006                   & \textbf{12.741 ± 0.001}  \\
                                  & 15                   & 8.526 ± 0.000  & 8.546 ± 0.000   & 8.554 ± 0.000    & 8.233 ± 0.001     & 1.220 ± 0.001             & -4.795 ± 0.001           & 12.688 ± 0.002                   & \textbf{12.742 ± 0.001}  \\
\multirow{3}{*}{SSIM (↑)}         & 5                    & 0.191 ± 0.005  & 0.131 ± 0.006   & 0.106 ± 0.004    & 0.111 ± 0.005     & 0.340 ± 0.000             & 0.338 ± 0.000            & 0.397 ± 0.001                    & \textbf{0.459 ± 0.002}   \\
                                  & 10                   & 0.261 ± 0.003  & 0.192 ± 0.003   & 0.167 ± 0.004    & 0.178 ± 0.003     & 0.340 ± 0.000             & 0.338 ± 0.000            & 0.425 ± 0.000                    & \textbf{0.462 ± 0.002}   \\
                                  & 15                   & 0.312 ± 0.001  & 0.256 ± 0.003   & 0.239 ± 0.003    & 0.256 ± 0.003     & 0.340 ± 0.000             & 0.338 ± 0.000            & 0.438 ± 0.000                    & \textbf{0.464 ± 0.001}   \\ \hline
\end{tabular}
\end{adjustbox}
\end{table}

To further quantitative evaluate the accuracy of recovered phase, we calculated the regression-phase SNR (\textit{rpSNR}) \cite{zhang2023retinex}, peak-SNR (\textit{PSNR}), and structure similarity (\textit{SSIM}) between the calculation and the ground truth as listed in Tab. \ref{Table 2}. For each group of SNR the experiments were repeated 10 times to obtain the mean value of the scores. According to Tab. \ref{Table 2}, the pd-qDPC gains better phase recovery quality among all three evaluation matrixes than that of Retinex TV-qDPC. The $L_2$-qDPC, Iso-qDPC, and TV-qDPC fail to reconstruct correct phase pattern (low \textit{rpSNR} and \textit{SSIM} scores) in this case.

\subsection{Quantitative phase targets}
\begin{figure}
    \centering
    \includegraphics[width = 0.95\textwidth,trim = 0 500 0 0,clip]{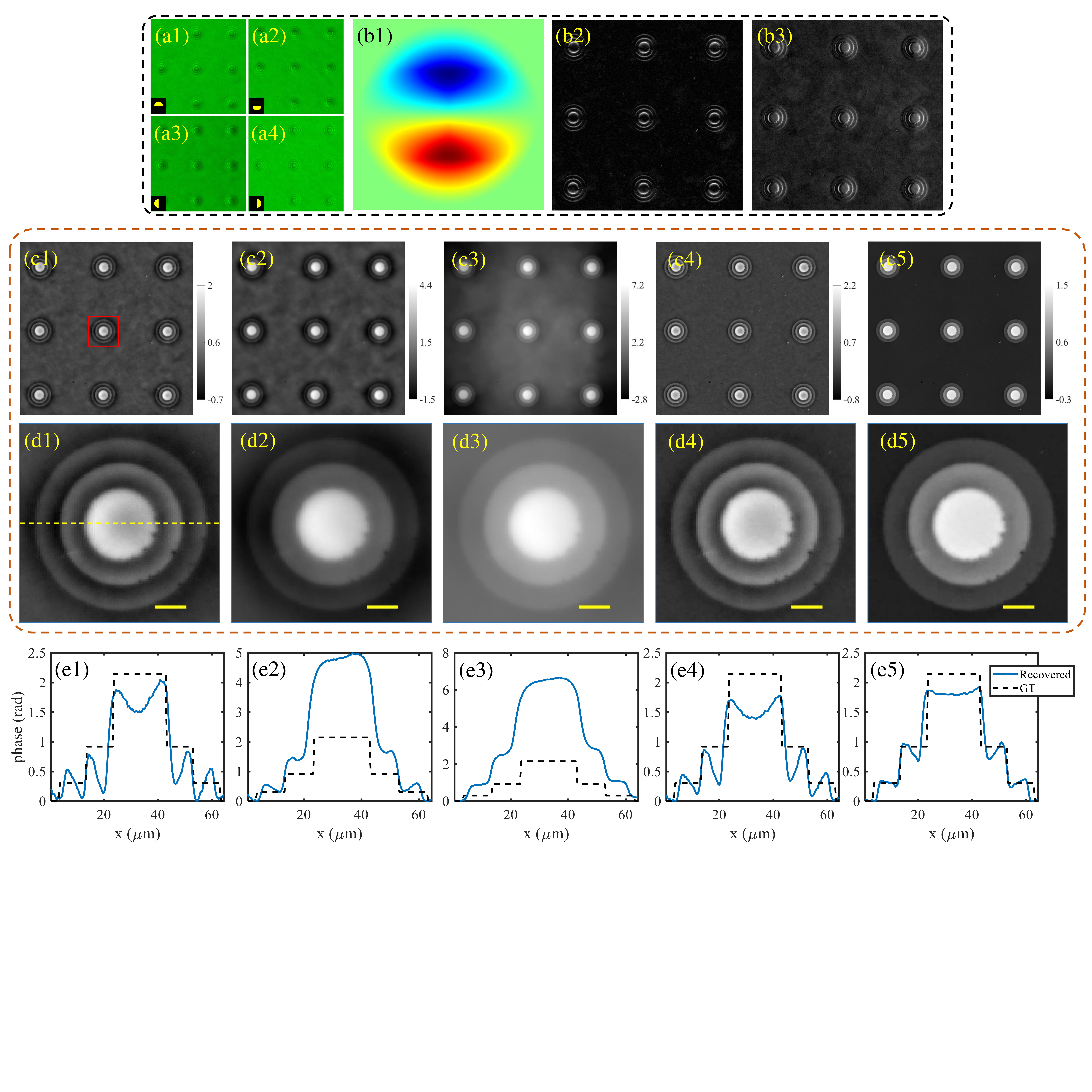}
    \caption{Experimental study on wedding-cake phase target. (a1) to (a2) are captured raw images. (b1) are PTF for up-bottom illumination. (b2) and (b3) are differential phase contrast images. (c1) to (c5) are reconstructued phase using $L_2\text{-qDPC}$, Iso-qDPC, TV-qDPC, Retinex TV-qDPC, and pd-qDPC, respectively. (d1) to (d5) are zoomed-in image in the red box. The scale-bar is 10 $\upmu\text{m}$ (e1) to (e5) are the quantitative phase profile along the dashed yellow line.}
    \label{Fig.6.}
\end{figure}

\begin{figure}
    \centering
    \includegraphics[width = 0.95\textwidth,trim = 0 310 0 0,clip]{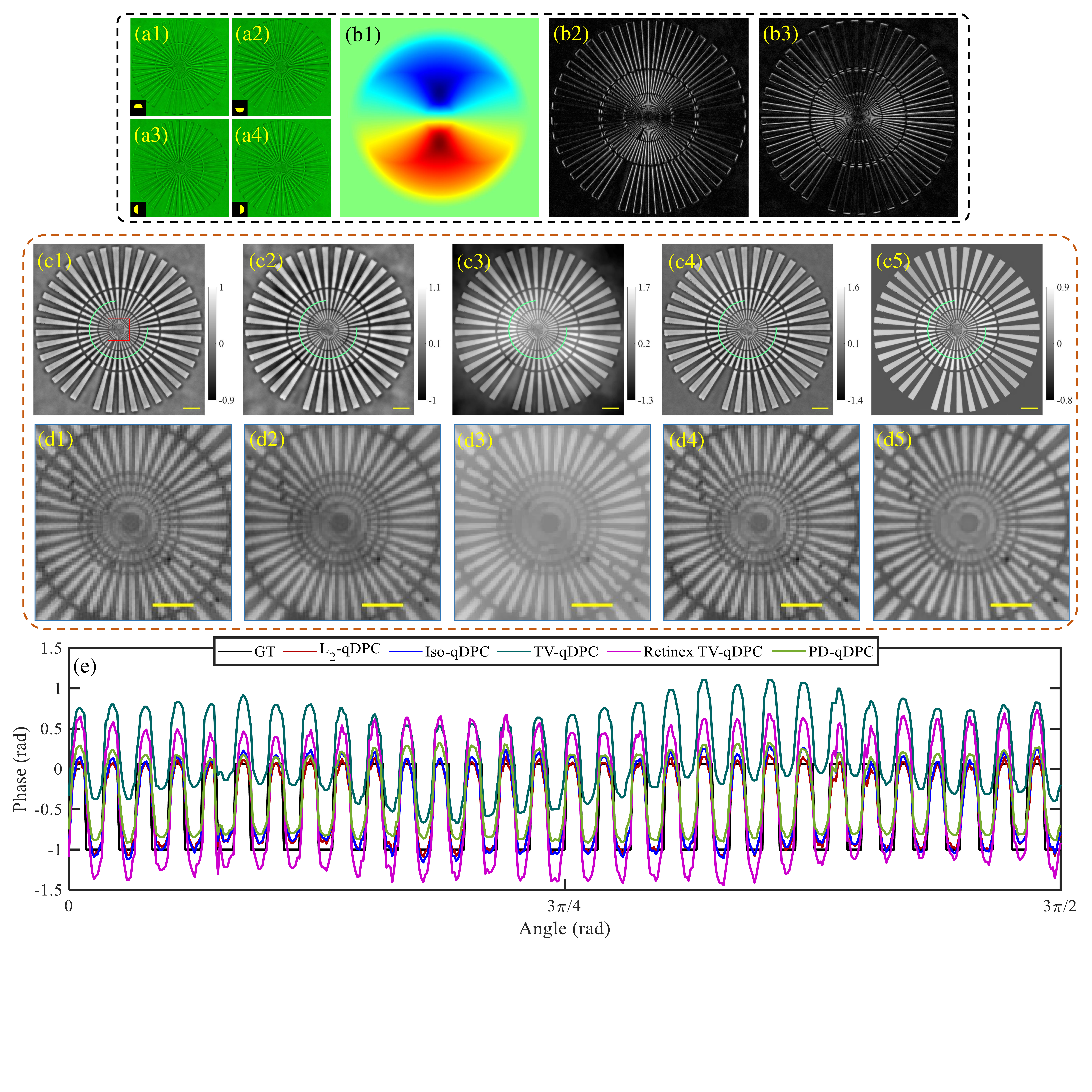}
    \caption{Experimental study on star-pattern phase target. (a1) to (a2) are captured raw images. (b1) are PTF for up-bottom illumination. (b2) and (b3) are differential phase contrast images. (c1) to (c5) are reconstructed phase using $L_2\text{-qDPC}$, Iso-qDPC, TV-qDPC, Retinex TV-qDPC, and pd-qDPC, respectively. The scale-bar is 40 $\mu m$ (d1) to (d5) are zoomed-in image in the red box. The scale-bar is 10 $\upmu\text{m}$ (e) is the quantitative phase profile along the green lines.}
    \label{Fig.7.}
\end{figure}

We conducted experimental study on two quantitative phase targets (QPT) including the wedding-cake pattern and the focal-star pattern, and shown in Fig. \ref{Fig.6.} and Fig. \ref{Fig.7.}, respectively. For the wedding-cake pattern, the captured oblique-illumination raw images are shown in Fig. \ref{Fig.6.} (a1) to (a4). The system parameters are: $\lambda=0.532 \upmu\text{m}$, NA = 0.3 and so does the illumination pupil. The magnification is $\times 10$, and the pixel size is $4 \upmu\text{m}$. The PTF for up-bottom illuminations are shown in Fig. \ref{Fig.6.} (b1). The DPC images are shown in Fig. \ref{Fig.6.} (b2) and (b3). In Fig. \ref{Fig.6.} (b3), one can see obverse mismatch of background intensity. 

Calculated phase patterns are shown in Fig. \ref{Fig.6.} (c1) to (c5), in which both Retinex TV-qDPC and pd-qDPC show background removal effect. Results of $L_2$-qDPC [Fig. \ref{Fig.6.} (c1)] and Iso-qDPC [Fig. \ref{Fig.6.} (c2)] are corrupted by background fluctuations and noise signals, while the 'white-cloud' effect is not pronounced as the as $L_2-\text{norm}$, known also as weight decay, in these two methods suppresses the value of the entire phase pattern. The TV-qDPC has significant 'white-cloud' effect as shown in Fig. \ref{Fig.6.} (c3) and (d3). 

Further as shown in the quantitative phase profile in Fig. \ref{Fig.6.} (e1) to (e5), our proposed pd-qDPC gains the best phase recovery quality in this group of study as the phase profile close to the GT value. 

A similar analysis was applied to Fig. \ref{Fig.7.} for the focal-star pattern. The system parameters were not changed, but we used a ring-shaped illumination pupil. The inner-radius is about 0.9 times the outer radius so that the PTF is presented as Fig. \ref{Fig.7.} (b1). Accordingly, our proposed pd-qDPC suppresses both the background fluctuations and noise signals as shown in Fig. \ref{Fig.7.} (c5) and (d5) without loss of valid structure information. The quantitative phase profile shows the high data fidelity of the pd-qDPC as shown in Fig. \ref{Fig.7.} (e). Traditional $L_2$-qDPC, Iso-qDPC, and TV-qDPC can suppress the noises but cannot handle the background mismatch, which is comment problem in DPC experiments. The Retinex TV-qDPC performs well in this case, while the TV regularization causes stair-case like artifact as shown in Fig. \ref{Fig.7.} (d4). 

\begin{figure}
    \centering
    \includegraphics[width = 0.95\textwidth,trim = 0 780 0 0,clip]{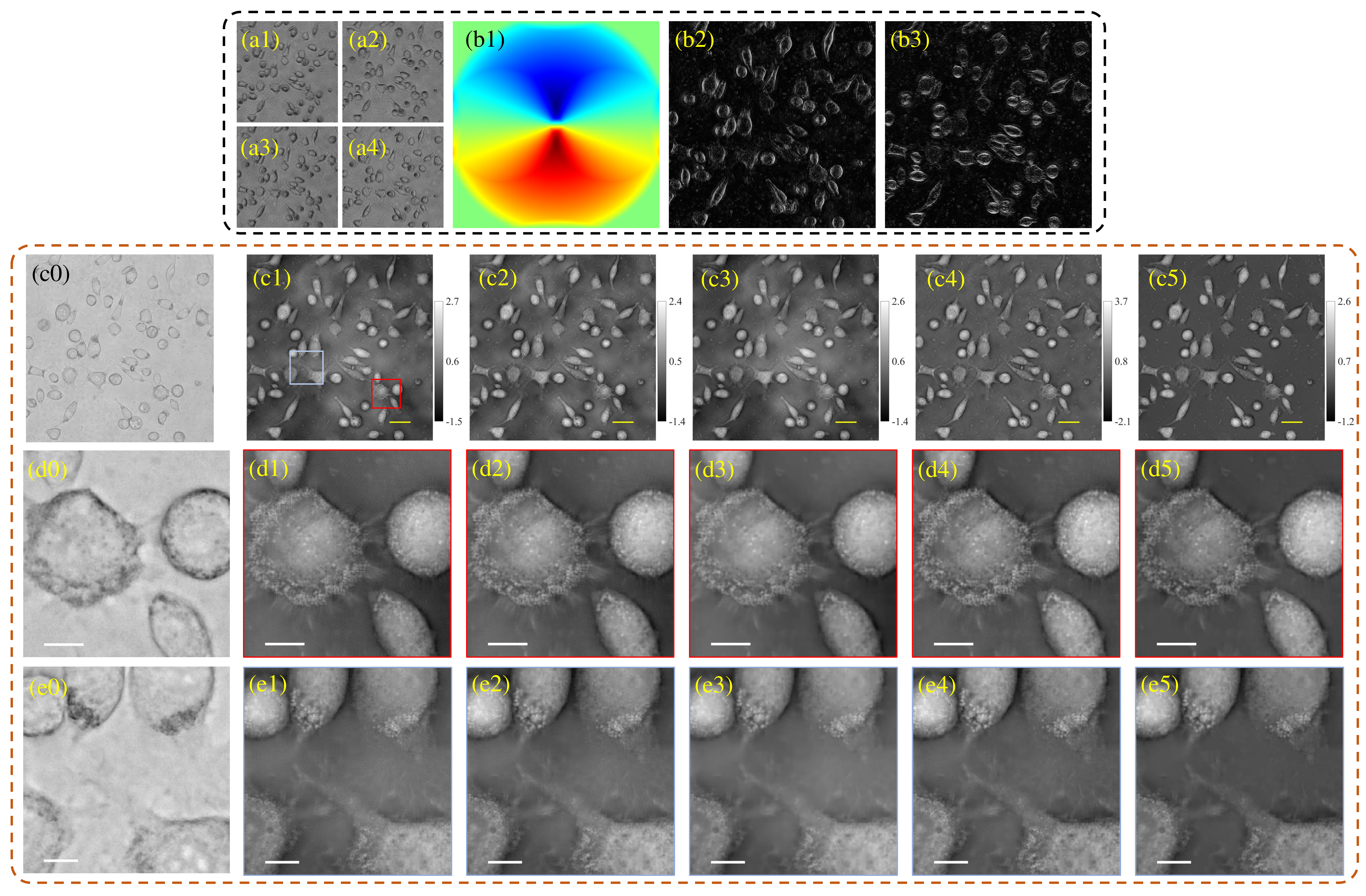}
    \caption{Experimental study on L929 cells cluster. (a1) to (a2) are captured raw images. (b1) are PTF for up-bottom illumination. (b2) and (b3) are differential phase contrast images. (c0), (d0), and (e0) are the bright field images. (c1) to (c5) are reconstructed phase using $L_2\text{-qDPC}$, Iso-qDPC, TV-qDPC, Retinex TV-qDPC, and pd-qDPC, respectively. The scale-bar is 40 $\mu m$ (d1) to (d5) are zoomed-in image in the red box. The scale-bar is 10 $\upmu\text{m}$ (e) is the quantitative phase profile along the green lines.}
    \label{Fig.8.}
\end{figure}

\subsection{L929 cells cluster}
We tested the pd-qDPC on L929 cells cluster as shown in Fig. \ref{Fig.8.}. In this experiment, NA = 0.5 and system magnification was $\times 20$. The pixel size of camera was 5.86$\upmu\text{m}$. We also used the ring-shaped illumination pupil. The maximum illumination NA = 0.5, and the inner radius was NA = 0.45, the PTF is shown in Fig. \ref{Fig.8.} (b1). 

The bright field images for L929 cells is shown in Fig. \ref{Fig.8.}(c0), and in (d0), (e0) for the zoomed-in area in the red and blue boxes. The phase images are shown in Fig. \ref{Fig.8.} (c1) to (c5), where the resolution is enhanced and many granular-membrane structures (lipid droplets maybe \cite{Osamu2021}) can be observed near the edge of the cells as shown in Fig. \ref{Fig.8.} (d1) to (d5). Comparing Fig. \ref{Fig.8.} (c5) to others, both qd-qDPC and Retinex TV-qDPC correct the background fluctuations among the entire field of view, while the qd-qDPC outperforms Retinex TV-qDPC since the background in Fig. \ref{Fig.8.} (c5) is more uniform than that of Fig. \ref{Fig.8.} (c4). Meanwhile, the noise signal is also suppressed, the thin and small membrane-structures are also maintained which can be seen from Fig. \ref{Fig.8.} (e5). 

In general, our proposed pd-qDPC outperforms State-of-the-art (SOTA) qDPC methods in background correction, noise suppression, and data fidelity. \textbf{Furthermore, as we will show in next section, the intermediate product of pd-qDPC is also of help for pattern recognition tasks including cell segmentation, cell contour tracking and edge detection}.
\section{Potential applications of pd-qDPC}\label{section 6}
\subsection{Automatic cell segmentation}
\begin{figure}
    \centering
    \includegraphics[width = 0.95\textwidth,trim = 0 1340 0 0,clip]{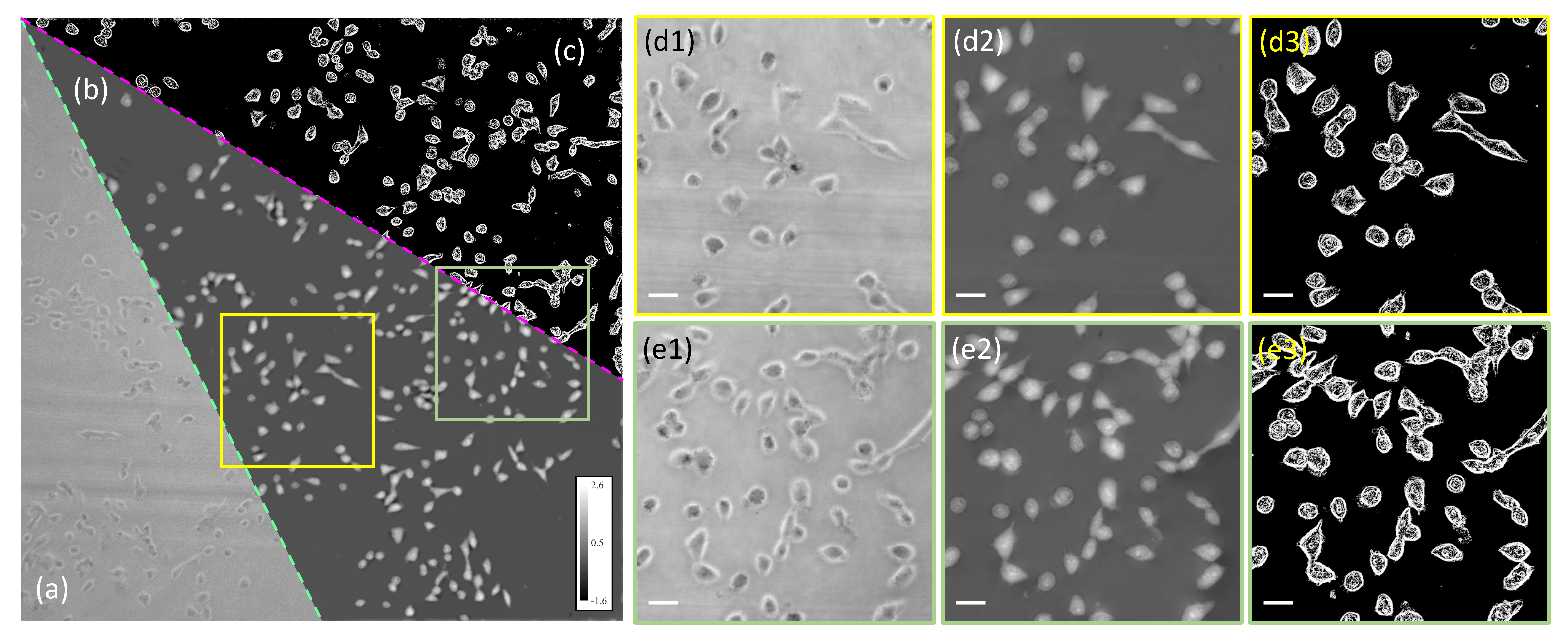}
    \caption{Experimental study on human gastric cancer cells (HCG-27) cluster. (a) bright field image. (b) quantitative phase image. (c) modulus of $\mathbf{\Psi}_n$. (d1) to (d3) are zoomed-in images for the yellow box area. (e1) to (e3) are zoomed-in images for the green box area. The scale bar is 20 $\mu m$}
    \label{Fig.9.}
\end{figure}

Our proposed pd-qDPC can automatically achieve cell contour segmentation during the phase reconstruction routine. The cell contour, or edge is directly related to the intermediate product of pd-qDPC given by the variables $\mathbf{\Psi}_n, (n=1,2, \cdots,  N)$, since the PTFs are naturally band-limited edge filter. Furthermore, the sparse regularization applied to the $\mathbf{\Psi}_n, (n=1,2, \cdots,  N)$ ensures a clear and black background pixels. For sparse distributed samples where the image contains most of background pixels, the cell contour segmentation effect of pd-qDPC becomes even obvious. 

To show this property, we conduct phase reconstruction of Hela cells cluster as shown in Fig. \ref{Fig.9.}. The quantitative phase image using pd-qDPC in Fig. \ref{Fig.9.} (b) shows its high-contrast phase image with even background. The modulus of the intermediate products $\mathbf{\Psi}_n$ is shown in Fig. \ref{Fig.9.} (c) where the cell contours are automatically segmented during the phase deconvolution. Moreover, the sparsity promotes the clear background where no noise signals and outliers appear. For cell segmentation purpose, we can further increase the penalty parameter $\alpha$ to achieve sparser segmentation results. The $\mathbf{\Psi}_n$ provides important geometrical information for cell shape than has potential application on unlabeled cell morphology analysis \cite{he2022morphology}. 

It is also worth noting that designing band-limited (noise-suppressed) edge filters can also be inspired by the PTF of DPC imaging, where edge filters can be generated using the convolution between a circular pattern and an odd function according to the formation of PTF. 

\subsection{Learning the PTF from edge response}

\begin{figure}
    \centering
    \includegraphics[width = 0.95\textwidth,trim = 0 1340 0 0,clip]{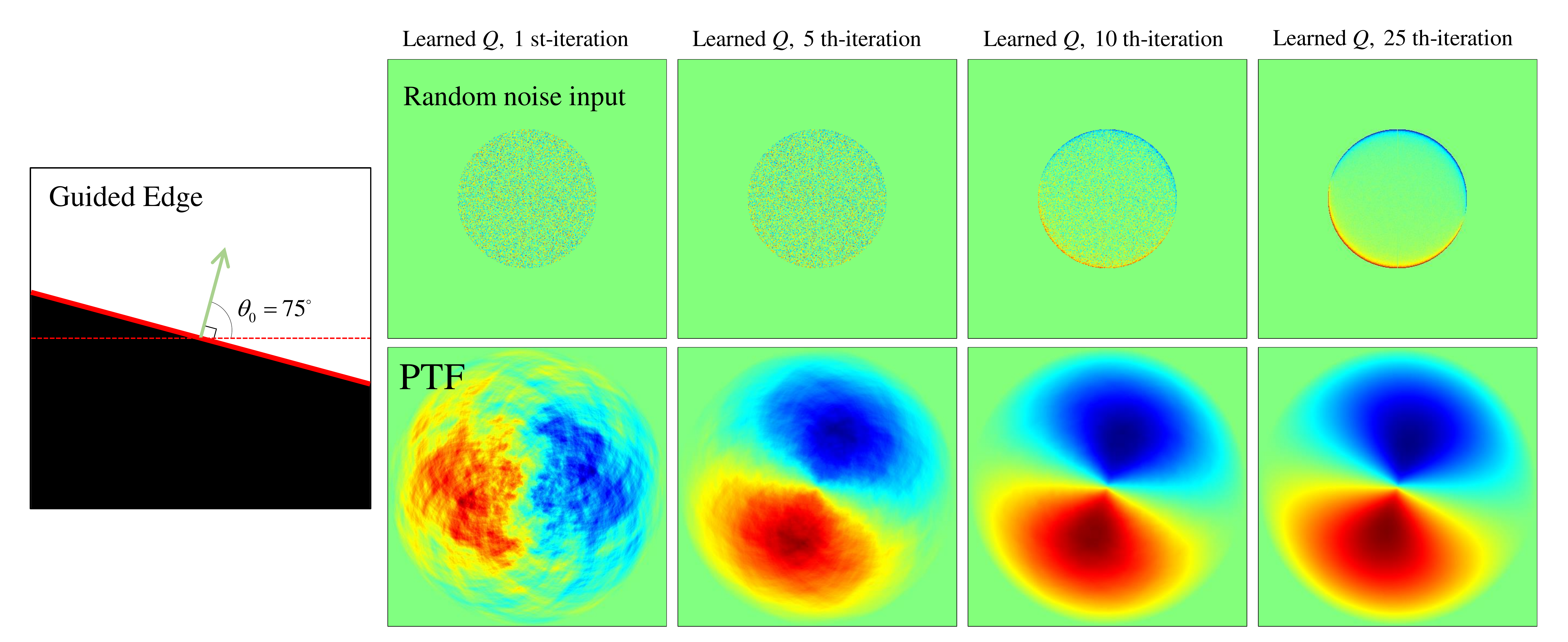}
    \caption{Learning the illumination pattern from the guided edge using edge-response cost function. The learned illumination pupils at 1, 5, 10 and 25 iterations are shown. The learning was performed within a normalized square area, while the radius of the pupil function is half of the width. \textbf{The trained illumination pupil automatically find the pupil-matching conditions without any given prior knowledge.}}
    \label{Fig.Learning.}
\end{figure}

pd-qDPC can guide the designing of PTF for DPC experiments. As described in Section \ref{section 3} the PTF is an edge filter, meaning that we can train a certain PTF by designing an edge-response strength cost function:
\begin{equation}
\begin{aligned}
    Q(\mathbf{f},\theta_0) & = \text{argmax}\left \| \left | \text{rot}\left ( \frac{1}{f_x},\theta_0  \right )\cdot  P(\mathbf{f} )\otimes\left (
P(\mathbf{f} )\cdot \left \{ Q(\mathbf{f})-\text{rot}\left [ Q(\mathbf{f}),\pi  \right ]  \right \}  \right )   \right |^2  \right \| ^2_2, \\
& \text{subject to} \left\| Q \right\|_{\infty} = 1.
\end{aligned}\label{(18)}
\end{equation}
The $\text{rot}(x,\theta)$ in Eq. (\ref{(18)}) denotes 2D rotating of function $x$ around the origin for $\theta$ angle, and $Q(\mathbf{f})-\text{rot}\left [ Q(\mathbf{f}),\pi  \right ]$  generates a odd-function w.r.t the origin. The absolute value and the $L_2$-norm measure the filter response of an edge $\text{rot}\left ( 1/f_x,\theta_0  \right )$ in the Fourier space according to the Plancherel's Theorem.

Since the pupil function $P$ is determined by the optical system, a certain illumination pattern can be learned by maximizing the cost function in Eq. (\ref{(18)}) using gradient descent method sharing the same ideal of training the convolutional neural network \cite{lecun1989backpropagation} and visualizing the convolution kernels \cite{zeiler2014visualizing}. 

For example, we trained a illumination pupil based on the guided-edge picture as shown in Fig. \ref{Fig.Learning.}. The input illumination pupil, $Q$, is a random Gaussian noise, while within a few iterations, $Q$ is convergent to a stable pattern as shown in the last column of Fig. \ref{Fig.Learning.}. This crescent-shaped illumination pupil was never reported in any publications related to qDPC. \textbf{Interestingly, the trained illumination pupil automatically finds the pupil-matching conditions \cite{fan2019optimal} where a strong circular appears along the edge of the illumination pupil without any given prior knowledge.}

In this way, we believe that the pd-qDPC brings together the pattern recognition and optical designing that can potentially benefit both fields.

\section{Discussions}\label{section 7}
\subsection{Ablations}
\begin{figure}
    \centering
    \includegraphics[width = 0.95\textwidth,trim = 0 1630 0 0,clip]{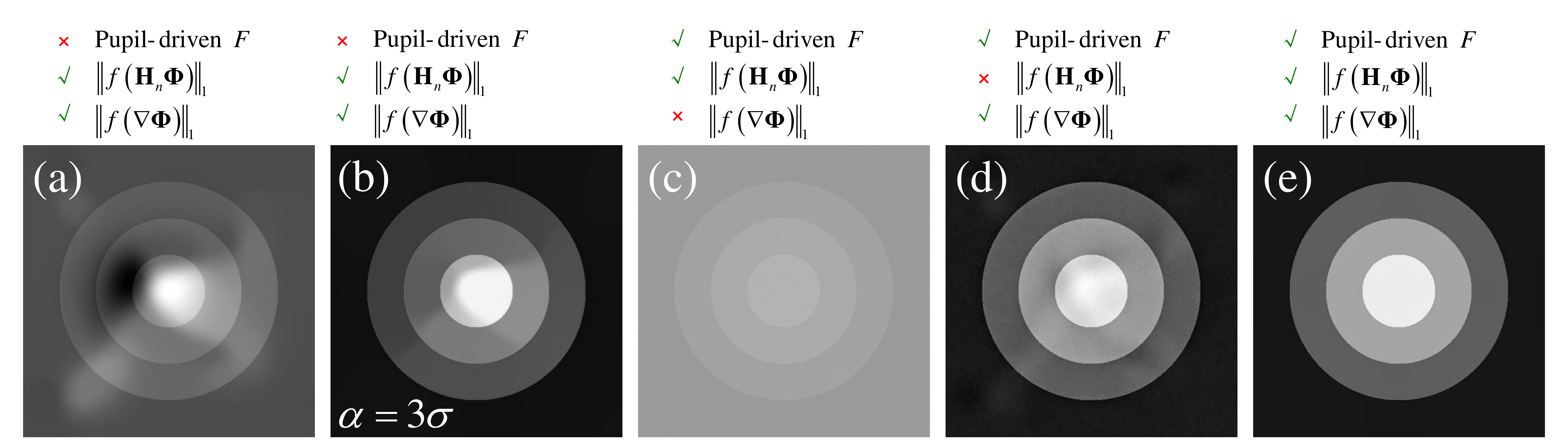}
    \caption{Ablation study for pd-qDPC. (a) to (e) are deconvolution results with different terms turned-off.}
    \label{Fig.10.}
\end{figure}

In this section, we performed an ablation study on the data fidelity term of pd-qDPC. Without the pupil driven term, the data fidelity term reduces to the traditional qDPC fidelity, without background-removal ability, as shown in Fig. \ref{Fig.10.} (a). Noise signal can be suppressed due to the gradient penalty term. In addition, by increasing the penalty parameter for the pupil-driven penalty term, the algorithm also shows its ability on background removal, as shown in Fig. \ref{Fig.10.} (b). On the other hand, the structure is over smoothed, which will further cause loss of tiny structure of phase images for smooth varying samples such as cells. 

\begin{table}
\renewcommand\arraystretch{1.5}
\caption{Comparison study on simulated data with background mismatch and different strength of noise.}
\label{Table 3}
\begin{adjustbox}{width=\columnwidth,center}
\begin{tabular}{lccccccc}
\hline
\multirow{3}{*}{Matrix} & \multicolumn{1}{l}{\multirow{3}{*}{SNR}} & Pupil driven fidelity           & \ding{53}    & \ding{53}                    & \checkmark & \checkmark & \checkmark \\
                        & \multicolumn{1}{l}{}                     & $f(\mathbf{H}_n \mathbf{\Phi})$ & \checkmark & \checkmark, large $\alpha$ & \checkmark & \ding{53}    & \checkmark \\
                        & \multicolumn{1}{l}{}                     & $f( \nabla\mathbf{\Phi})$       & \checkmark & \checkmark                 & \ding{53}    & \checkmark & \checkmark \\ \hline
\multirow{3}{*}{rpSNR}  & 5                                        & \multirow{9}{*}{-}              & 1.702 ± 0.022             & 11.484 ± 0.046                            & 4.985 ± 0.072             & 19.701 ± 0.171            & \textbf{20.423 ± 0.089}   \\
                        & 10                                       &                                 & -1.603 ± 0.015            & 6.903 ± 0.026                             & 6.819 ± 0.045             & 16.482 ± 0.096            & \textbf{17.645 ± 0.062}   \\
                        & 15                                       &                                 & -3.463 ± 0.007            & 1.438 ± 0.014                             & 8.780 ± 0.035             & 13.281 ± 0.016            & \textbf{17.070 ± 0.044}   \\
\multirow{3}{*}{PSNR}   & 5                                        &                                 & 8.315 ± 0.015             & 12.073 ± 0.007                            & 10.141 ± 0.034            & \textbf{12.714 ± 0.005}   & 12.334 ± 0.010            \\
                        & 10                                       &                                 & 5.914 ± 0.012             & 10.925 ± 0.009                            & 10.895 ± 0.016            & 12.582 ± 0.005            & \textbf{12.641 ± 0.003}   \\
                        & 15                                       &                                 & 4.373 ± 0.006             & 8.143 ± 0.010                             & 11.506 ± 0.009            & 12.318 ± 0.002            & \textbf{12.614 ± 0.002}   \\
\multirow{3}{*}{SSIM}   & 5                                        &                                 & 0.368 ± 0.000             & 0.433 ± 0.001                             & 0.203 ± 0.001             & 0.439 ± 0.002             & \textbf{0.453 ± 0.002}    \\
                        & 10                                       &                                 & 0.355 ± 0.000             & 0.374 ± 0.001                             & 0.211 ± 0.001             & 0.444 ± 0.001             & \textbf{0.457 ± 0.001}    \\
                        & 15                                       &                                 & 0.347 ± 0.000             & 0.368 ± 0.000                             & 0.224 ± 0.000             & 0.440 ± 0.000             & \textbf{0.455 ± 0.001}    \\ \hline
\end{tabular}
\end{adjustbox}
\end{table}

If we further turned off the gradient penalty term as shown in Fig. \ref{Fig.10.} (c), the deconvoluted image was severely degraded by the present of noise since the kernel is ill-conditioned. The gradient penalty suppresses noise signal by minimizing the variation, similar to the TV regularization. 

With the pupil-driven fidelity and pupil-driven penalty turned off, our pd-qDPC works similar to the Retinex TV-qDPC where the background fluctuation can be efficient suppressed as shown in Fig. \ref{Fig.10.} (d). However, the gradient penalty only cannot ensure the sparsity of gradient of DPC images as analysed in Section \ref{section 3}, as such, a "white-cloud" effect also appears to the phase pattern. By turning on the pupil-driven penalty, we obtained high quality full-field phase deconvolution without the impact of noise signal and a background intensity. Statistical analysis based on simulation data is listed in Tab. \ref{Table 3}, showing that each term in the cost function is efficient and work properly.

Especially, the pupil-driven penalty, $ \left \|f(\mathbf{H_n \mathbf{\Phi}}) \right \|_1$, is effective because it has lower energy for idea DPC images than for noise-corrupted DPC images as stated in Section \ref{section 3}. It favors the idea DPC images as the idea DPC image is so sparser than the noise-corrupted DPC image that $ \left \|f(\mathbf{H_n \mathbf{\Phi}}) \right \|_1$ can simply distinguish the idea DPC image from that of noise-corrupted DPC image. In this manner, the pupil-driven penalty strongly rectifies the qDPC deconvolution making the it well-conditioned.

\subsection{Parameters tuning}
The pd-qDPC cost function contains three parameters, $\omega$, $\alpha$, and $\beta$. The value of $\omega$ controls the RST operator and has little impact to the deconvolution results according to Fig. \ref{Fig.4.} (b). $\omega = 10$ is suitable for most of qDPC applications based on our simulation and experimental studies. Parameters $\alpha$ and $\beta$ can be automatically determined using Eq. \ref{(17)}. The adaptive noise sensor estimates the noise strength beyond the band-limit of PTF and thus can correctly evaluate the noise strength even for image without noise, the noise sensor also provide reliable guidance for fine tuning the model parameters, and fortunately, tuning these two is not time-consuming and is way easier than tuning the layout of optical system as the execution speed of pd-qDPC is fast, which costs less than 10 seconds for input images with 2048 by 2048 pixels. (Intel(R) Xeon(R) W-10855M CPU @ 2.80GHz 2.81 GHz, 32GB Memory)
\section{Concluding remarks}\label{section 8}

In conclusion, we proposed here a new qDPC reconstruction framework, termed pupil-driven qDPC (pd-qDPC), that embeds the inborn edge filters of DPC-the phase transmission function, to facilitate the qDPC's phase deconvolution ability for robust and high-fidelity phase retrieval. 

We proved that the PTF is an edge filter, and formulated a non-convex model for the inverse problem to achieve automatic background removal, noise suppression, and cell edge segmentation, simultaneously. Both simulation and experimental results show the superiority of the pd-qDPC in terms of phase reconstruction quality and implementation efficiency compared with existing techniques. Importantly, the PTF for qDPC can be also learned from the edge-response cost function that potentially benefits the qDPC applications.

Since pd-qDPC is purely mathematical, it doesn’t require any additional physical prior information and modification to the optical system, in contrast with existing methods. No background calibration/removal is needed for sample induced background, and a good phase pattern can be reconstructed even if the DPC data is degenerated by the severe noise. pd-qDPC is of strong interpretability and generalization, that can be used to boost all qDPC experiments including but not limited to 3D qDPC image \cite{chen20163d,ling_BOE} and single-shot qDPC image.

\subsubsection{Acknowledgements}
\
\newline
This research was funded by Universiteitsfonds Limburg SWOL 2022 (CoBes22.041); Research Fund of Anhui Institute of Translational Medicine (2021zhyx-B16); Key Research and Development Program of Anhui Province (2022a05020028); China Scholarship Council (CSC) (201908340078).

\subsubsection{Declaration of interests}
\
\newline
The authors declare that they have no known competing financial interests or personal relationships that could have appeared to influence the work reported in this paper.
%
%
%
%

\printbibliography
\end{document}